\documentclass[a4paper,12pt]{amsart}
\usepackage{amsmath,amssymb,latexsym,amsfonts,amscd}

\title{Small bound for birational automorphism groups of algebraic varieties \\
\rm{(with an Appendix by Yujiro Kawamata)}}
\author{De-Qi Zhang}

\thanks{The author is supported by an Academic Research Fund of
NUS
\newline
2000 Mathematics Subject Classification 14J50, 14E07}
\newcommand{\CC}{{\mathbb C}}

\newcommand{\PP}{{\mathbb P}}

\newcommand\Alb{\text{\rm Alb}}
\newcommand\alb{\text{\rm alb}}

\newcommand\OO{{\mathcal{O}}}

\newcommand\ZZ{{\mathbb{Z}}}

\newcommand\Sing{{\text{\rm Sing}}}
\newcommand\Imm{{\text{\rm Im}}}
\newcommand\id{{\text{\rm id}}}
\newcommand\diag{{\text{\rm diag}}}
\newcommand\ord{{\text{\rm ord}}}
\newcommand\Tr{{\text{\rm Tr}}}
\newcommand\Aut{{\text{\rm Aut}}}
\newcommand\Bir{{\text{\rm Bir}}}
\newcommand\variety{{\text{\rm variety}}}
\newcommand\group{{\text{\rm group}}}

\newcommand\GL{{\text{\rm GL}}}
\newcommand\SL{{\text{\rm SL}}}
\newcommand\lcm{{\text{\rm lcm}}}
\newcommand\torsion{{\text{\rm torsion}}}
\newcommand\Ker{{\text{\rm Ker}}}
\newcommand\cartier{{\text{\rm cartier}}}

\newtheorem{thm}{Theorem}[section]
\newtheorem{lem}[thm]{Lemma}
\newtheorem{cor}[thm]{Corollary}
\newtheorem{prop}[thm]{Proposition}
\newtheorem{claim}[thm]{Claim}
\theoremstyle{definition}

\newtheorem{setup}[thm]{}

\newtheorem{conj}[thm]{Conjecture}
\newtheorem{rem}[thm]{Remark}
\theoremstyle{remark}

\begin{document}
\begin{abstract}
We give an effective upper bound of $|\Bir(X)|$ for the birational
automorphism group of an irregular $n$-fold (with $n = 3$) of
general type in terms of the volume $V = V(X)$ under an ''albanese
smoothness and simplicity'' condition. To be precise, $|\Bir(X)| \le
d_3 V^{10}$. An optimum linear bound $|\Bir(X)| \le \frac{1}{3}
\times 42^3 V$ is obtained for those $3$-folds with non-maximal
albanese dimension. For all $n \ge 3$, a bound $|\Bir(X)| \le d_n
V^{10}$ is obtained when $\alb_X$ is generically finite, $\alb(X)$
is smooth and $\Alb(X)$ is simple.
\end{abstract}
\maketitle
\pagestyle{myheadings} \markboth{\hfill D. -Q.
Zhang\hfill}{\hfill Automorphism groups of algebraic varieties\hfill}
%
%
%
%
\section{\bf Introduction}
We work over the field $\CC$ of complex numbers. Let $X$ be a
normal projective $n$-fold of general type. $X$ is {\it minimal}
if the canonical divisor $K_X$ is nef and $X$ has at worst
terminal singularities (see Kawamata-Matsuda-Matsuki \cite{KMM},
Kollar-Mori \cite{KM}). It is known that $|\Aut(X)| \le 42
\deg(K_X)$ when $\dim X = 1$ (Hurwitz), and $|\Aut(X)| \le (42
K_X)^2$ when $X$ is a minimal surface of general type (Xiao
\cite{Xi94}, \cite{Xi95}). See also works of Andreotti,
Howard-Sommese \cite{HSo}, Huckleberry - Sauer \cite{HS}, Corti
\cite{Co} and probably many others.

For Gorenstein minimal $n$-folds $X$ of general type, Szabo
\cite{Sz} gave a polynomial upper bound of $|\Bir(X)|$, though its
degree is quite huge. See also works of Catanese - Schneider
\cite{CS}, Xiao \cite{Xi96}, Cai \cite{Cai} and Kovacs \cite{Kov}.

In this paper, we try to get a more realistic upper bound (of order)
for the whole birational automorphim group $\Bir(X)$. If fact, we
consider irregular varieties $X$ and, {\it without} assuming the
minimality of $X$, obtain a linear (resp. degree $10$) bound in the
volume $V(X)$ when $n = 3$ (resp. for all $n \ge 3$) if the albanese
map is not generically finite onto a variety of general type (resp.
is generically finite onto a smooth variety contained in the simple
Albanese variety); for details, see below.

Our approach here of considering the albanese morphism $\alb_X : X
\rightarrow \Alb(X)$ is sort of the generalization of the idea as in
Xiao \cite{Xi90} and Catanese-Schneider \cite{CS} to find a
$G$-equivariant pencil. The universal property of $\Alb(X)$
guarantees that any action of a group $G$ on $X$ induces a canonical
action on $\Alb(X)$ such that $\alb_X$ is $G$-equivariant though $G$
might not act faithfully on the latter. When $\alb_X$ is not
generically finite onto a variety of general type, Ueno's result
\cite{Ue} says that $Y: = \alb_X(X)/B_0$ is of general type, where
$B_0$ is the identity connected component of the subtorus of
$\Alb(X)$ stabilizing $\alb(X)$ with respect to its translation
action. So one can apply the weak positivity of the relative
dualizing sheaf due to Fujita \cite{Fuj}, Kawamata \cite{Ka82b} and
Viehweg \cite{Vi82b}, then use the volume $V(X)$ to give an optimum
upper bound of $V(Y) V(F)$ and finally reduce to the cases of the
fibre and base ($F$ being the general fibre of $X \rightarrow Y$).

The ''best'' case where $\alb_X$ is generically finite onto $W =
\alb_X(X)$ of general type, turns out to be the hardest one. It is
supposed to be the ''best'' situation because then $|6K_X|$ is
birational according to Chen - Hacon \cite{CH} Corollary 5.3 . It
may then improve the coefficient in Szabo's bound, but not the power
of $V = V(X)$.

Therefore, we have to explore closely good properties of the abelian
variety $\Alb(X)$ as a torus. Certainly, after quotient away the
translations, the image $\overline{G}$ of $G = \Aut(X) \rightarrow
\Aut_{\variety}(\Alb(X))$ can be thought of as a subgroup of
$\GL_{2q}(\ZZ)$ with $q = q(X) = \dim \Alb(X)$, via the rational
representation on the first integral homology of $\Alb(X)$. But
Feit's upper bound for the order of a finite subgroup in
$\GL_{2q}(\ZZ)$ with $2q > 10$, is a huge number $(2q)! \, 2^{2q}$
(see Friedland \cite{Fr}) and it is even attainable by the
orthogonal group $O_{2q}(\ZZ)$. Or one can apply Jordan's Lemma
\ref{Jordan} below, together with Xiao's linear bound for abelian
subgroups, but then the constant $J_{2q}$ depends on $q$ and hence
on the volume of $X$, rather than on the dimension of $X$.

So we have to catch the missing information of $X$ when passing to
$\Alb(X)$. To do so, we first give a linear bound of Betti numbers
of $X$ in terms of the volume of $X$ in Theorem \ref{Bbetti} below,
then bound the exponent $\exp(G)$ using the classical Lefschetz
fixed point formula in Lemma \ref{Lefix} below, and finally bound
$|G|$ itself. In the process, one needs a technical finiteness
condition on the fixed locus $\Alb(X)^{\overline{G}}$ as in Theorem
\ref{Th2}. We remark that this condition is automatically satisfied
when $\Alb(X)$ is simple. Hopefully, one will be able to remove this
restriction and eventually have the desired small bound for all
irregular varieties.

One reason of our considering irregular varieties is that these ones
are expected to have a bigger $\Aut(X)$ compared with the general
ones.

The advantage of using the albanese map for irregular varieties is
that one could bound $\Aut(X)$ inductively by reducing either to
the cases of the fibre and base of a fibration, or to the case
where $\alb_X$ is generically finite onto a variety of general
type (which is covered by Theorem \ref{Th2} to some extent), since
$\alb_X$ always carries an action of $G$ on $X$ over to actions on
the base and fibre of $\alb_X$. This method is applicable in all
dimensions. See Theorem \ref{Thconj} below for details.

We now state the main results of the paper. When $\alb_X$ is not
generically finite onto a variety of general type, one obtains the
optimum linear bound; see the Remark in the Appendix; see
\ref{setup1}.

\begin{thm}\label{Th1}
Let $X$ be a smooth projective $3$-fold of general type
with irregularity $q(X) \ge 4$. Let $V = V(X)$ be the volume of $X$
and $G = \Bir(X)$ ($\ge \Aut(X)$)
the birational automorphism group of $X$.

\begin{list}{}{
\setlength{\leftmargin}{10pt}
\setlength{\labelwidth}{6pt}
}
\item [\rm{(1)}]
Suppose that $X$ is not of maximal albanese dimension,
i.e., $\dim$ $\alb_X(X) < \dim X$.
Then
$|G| \le \frac{1}{3} \times 42^3 \, V$.

\item [\rm{(2)}]
More generally, suppose only that $\alb_X : X \rightarrow \Alb(X)$
is not generically finite
onto a $3$-fold of general type in $\Alb(X)$,
i.e., the Kodaira dimension $\kappa(\alb_X(X)) < \dim X$.
Then
$|G| \le \frac{1}{3} \times 42^3 \, V$.
\end{list}
\end{thm}

When $\alb_X$ is generically finite onto a variety of general type,
i.e., $\kappa(\alb_X(X)) = \dim X$, one has the following slightly
bigger bound, but the degree of the bound is still a constant
independent of $\dim X$.

\begin{thm}\label{Th2}
Let $X$ be a smooth projective $n$-fold ($n \ge 3$) of general type.
Let $V = V(X)$ be the volume of $X$ and $G = \Bir(X)$.
Assume the following conditions:

\begin{list}{}{
\setlength{\leftmargin}{10pt}
\setlength{\labelwidth}{6pt}
}
\item [\rm{(i)}]
$\alb_X$ is generically finite onto $W$ ($= \alb(X)$) of general type,

\item [\rm{(ii)}]
$W$ is smooth, and

\item [\rm{(iii)}]
either $A: = \Alb(X)$ is a simple abelian variety; or $G$ induces an
action on $A$ such that the fixed locus $A^g$ for every $\id \ne g
\in G$ is a non-empty finite set unless $g$ is a translation of $A$.
\end{list}
\noindent
Then there is a constant $d_n$ (indepedent of $X$)
such that $|G| \le d_n V^{10}$.
\end{thm}

Combining Theorems \ref{Th1} and \ref{Th2}, one obtains a reasonably
small upper bound of $|\Bir(X)|$ for some irregular $3$-folds of
general type.

One ingredient for the proofs of results above is the following
linear bound of Betti numbers in terms of the volume.

This might give another simple proof of Xiao's linear bound of
$\ord(g)$ in terms of $V(X)$ for every $g \in \Aut(X)$ with
$\ord(g)$ prime under the assumption that $K_X$ is ample; see Remark
\ref{remBb}.

\begin{thm}\label{Bbetti}
Let $X$ be a smooth projective $n$-fold with ample $K_X$. Then there
is a constant $a_n$, independent of $X$, such that the Betti numbers
$B_i(X) \le a_n K_X^n$ for all $i$.
\end{thm}

\begin{rem}\label{remIntro}
\noindent (1) The constants $d_n$ and $a_n$ in the results above
are computable from the proof. The main contributors towards them
are some of the existing constants: $x_n$ in Xiao's Theorem
\ref{XiaoTh}, $J_n$ the Jordan constant in Lemma \ref{Jordan},
$r_n$ of Angehrn - Siu in Lemma \ref{smalldiminv} and $h_n$ of
Heier in Lemma \ref{eX}.

\noindent (2) Note that the degree ($= 10$) of the polynoimal in
$V(X)$ in Theorem \ref{Th2}, is independent of the dimension $n$
of $X$. We propose the following, which is somewhat more general
than the one in Xiao \cite{Xi94} for smooth and minimal $X$ and
for $\Aut(X)$:
\end{rem}

\begin{conj}\label{conj}
There is a constant $\delta_n$ such that the following holds for
every smooth projective $n$-fold of general type, where $V = V(X)$
is the volume:
$$|\Bir(X)| \le \delta_n V.$$
\end{conj}

When $n \le 2$, Conjecture \ref{conj} is confirmed by Hurwitz, Xiao
\cite{Xi94} and Xiao \cite{Xi95}, and one can take $\delta_n =
(42)^n$. Theorem \ref{Th1} confirms the conjecture for $n = 3$ with
some additional assumption. The result below is a further evidence
for arbitrary dimension.

\begin{thm}\label{Thconj}
Let $X$ and $Y$ be smooth projective varieties of general type and
let $f : X \rightarrow Y$ be a surjective morphism with connected
general fibre $F$. Suppose that $G = \Bir(X)$ acts regularly and
faithfully on $X$ and $G$ acts on $Y$ so that $f$ is
$G$-equivariant.
\par
If Conjecture \ref{conj} holds for both $Y$ and $F$, then it also
holds for $X$.
\end{thm}

\begin{setup} {\bf Terminology and Notation.} \label{setup1}
\par
For a normal projective variety $X$ of dimension $n$, we say that
$X$ is {\it minimal} if the canonical divisor $K_X$  is nef and if
$X$ has at worst terminal singularities; see \cite{KMM}, \cite{KM}.
\par
For a divisor $D$ on $X$, we define the {\it volume} of $D$ as $V(D)
= \lim \sup_{s \rightarrow \infty} h^0(X, s D)/(s^n /n!)$; see
\cite{La}. $V(X) := V(K_X)$ is called the {\it volume} of $X$, a
birational invariant. Note that $V(X) = K_X^n$ when $X$ is minimal.

Let $X$ be a variety and $G \le \Aut(X)$. Let $G_x = \{g \in G \, |
\, g(x) = x\}$ be the {\it stabilizer} subgroup. For $H \le G$,
define the {\it fixed locus} $X^H := \{x \in X \, | \, h(x) = x$ for
some $\id \ne h \in H\}$. For $g \in G$, define the {\it fixed
locus} $X^g := \{x \in X \, | \, g(x) = x\}$. Thus normally,
$X^{\langle g \rangle} \supseteq X^g$.

For an abelian variety $A$, denote by $\Aut_{\variety}(A)$ the set
of automorphisms of $A$ as a variety, and $\Aut_{\group}(A)$ the
subgroup of bijective homomorphisms.

When $H \le G$ and $g \in G$, denote by
$C_H(g) = \{h \in H \, | \, gh = hg\}$
the {\it centralizer}. Define the {\it exponent}
$\exp(G) = \lcm\{\ord(g) \, | \, g \in G\}$.
\end{setup}

\noindent {\bf Acknowledgment.} This work was started when I was
visiting National Taiwan University in the summer of the year 2005.
I am very grateful to Professor Alfred Chen for the discussion at
the initial stage, Professor Kawamata for kindly writing an Appendix
with which the old bound of degree 2 in the old version of Theorem
\ref{Th1} has been optimised to the current one, Professor Oguiso
for the reference \cite{Ue2}, and the referee for going the extra
mile to help me in removing the ambiguous parts and better the
presentation of the paper.
\section{\bf Preliminary results}
We first show how $K_X^n$ controls (or is controlled by) other
invariants linearly. For its generalization to arbitrary $n$, see
Proposition \ref{qp}. We remark that in the assertion (4) below, one
can take $r = 2 + n(n+1)/2$ or bigger according to Angehrn-Siu
\cite{ASiu}, or Kollar \cite{Ko97} Theorem $5.8$.

\begin{lem} \label{smalldiminv}
Let $X$ be a Gorenstein minimal projective $n$-fold ($n \ge 2$)
of general type
and $f : X' \rightarrow X$ a resolution.

\begin{list}{}{
\setlength{\leftmargin}{10pt}
\setlength{\labelwidth}{6pt}
}
\item[\rm{(1)}]
Suppose $n = 2$. Then $q(X) \le p_g(X) \le \frac{1}{2}K_X^2 + 2$,
and $K_X^2 \le 9(1 + p_g(X))$.

\item[\rm{(2)}]
Suppose $n = 3$. Then $p_g(X) \le K_X^3 + 2$ and $q(X) \le 2 + 50
K_X^3$. Also $\chi(\OO(K_X)) = f^* K_X . c_2(X')/24 \ge K_X^3/72 >
0$ and
$$P_m(X) = \frac{1}{12}m(m-1)(2m-1) K_X^3 + (2m-1) \chi(\OO(K_X))
\,\,\,\, (m \ge 2).$$
\item[\rm{(3)}]
Suppose $n = 3$ and $p_g(X) > 0$. Then
$K_X^3 \le 144 p_g(X)$.

\item[\rm{(4)}]
Suppose $n \ge 2$. Then for every $r \ge 4$ such that $|rK_X|$ is
base point free and $|(r+1) K_X|$ is non-empty we have $P_{r}(X) \le
n + r^{n}K_X^n/2$.
\end{list}
\end{lem}

\begin{proof}
(1) Now $X$ is smooth for ''terminal'' means smooth in dimension 2.
Note that $0 < \chi(\OO_X) = 1 - q(X) + p_g(X)$ for surfaces of
general type. Then the first two inequalities follow from the
Noether inequality. The last one follows from the Miyaoka-Yau
inequality in \cite{Mi} Theorem 1.1 and \cite{Yau}: $K_X^2 \le
3c_2(X)$ and the calculation: $1 + p_g(X) \ge \chi(\OO_X) = (K_X^2 +
c_2(X))/12 \ge K_X^2/9$.

(2) By Chen \cite{C} Theorem 3, we have $p_g(X) \le 2 + K_X^3$. By
Lee \cite{Lee00}, $|4K_X|$ is base point free and we take its
general memeber $S$ which is smooth (noting that $\Sing(X)$ is
finite because $X$ is terminal). Consider the exact sequence
$$0 \rightarrow \OO_X(-S) \rightarrow \OO_X
\rightarrow \OO_S \rightarrow 0.$$
Note that $H^i(X, \OO_X(-S)) = 0$ for $i = 1, 2$ by
Kawamata-Viehweg vanishing. Taking cohomology of the exact
sequence, we get $q(X) = q(S)$. Now by (1), $q(S) \le p_g(S) \le
\frac{1}{2}K_S^2 + 2$. Substituting in $K_S = (K_X + S)|S = 5K_X |
S$ and $K_S^2 = 100 K_X^3$, we get the first part of (2). The
second part of (2) follows from the Riemann-Roch formula in Reid
\cite{YPG} and Miyaoka-Yau inequality.

(3) By the proof of Hacon \cite{Hac} page 6, under the assumption
that $p_g(X) > 0$, one has $\chi(\omega_X) \le 2 p_g(X)$. Now (3)
follows from (2).

(4) Consider first the case $\dim X = 2$. The system $|rK_X|$ is
base point free for all $r \ge 4$, by Bombieri \cite{Bo}. By (1),
$\chi(\OO_X) \le 1 + p_g(X) \le 3 + K_X^2/2$. The pluri-genus
formula says that $P_r(X) = r(r-1)K_X^2/2 + \chi(\OO_X) \le 3 +
\frac{1}{2} (1 + r(r-1)) K_X^2$. One can verify that (4) is true.

Now assume that $n \ge 3$. Let $X_{n-1}$ be a general member of
$|rK_X|$. By Kollar-Mori \cite{KM} Lemma 5.17, $X_{n-1}$ is again
terminal with $K_{X_{n-1}} = (1 + r) K_X |{X_{n-1}}$. Let $X_i$ be
the the intersection  of $n-i$ general members of $|rK_X|$.
Inductively, we see that each $X_i$ is terminal with nef and big
$K_{X_i} = (1 + (n-i)r) K_X|X_{i}$. Also $C = X_1$ is a connected
smooth curve; note that $C$ is linearly equivalent to the nef and
big divisor $rK_X | X_2$ on the terminal, which is smooth in
dimension 2, surface $X_2$; one may also refer to Hartshorne
\cite{Hart} Ch III, Exercise 11.3, for the connectedness of each
$X_i$. Inductively, we obtain $h^0(X, rK_X) \le 1 + h^0(X_{n-1},
rK_X |X_{n-1}) \le \dots \le (n-1) + h^0(C, rK_X | C)$, by
considering the exact sequence of cohomologies coming from
$$0 \rightarrow \OO_{X_{i+1}}
\rightarrow \OO_{X_{i+1}}(rK_X | X_{i+1}) \rightarrow
\OO_{X_i}(rK_X|X_i) \rightarrow 0.$$ The divisor $D = rK_X | C$ on
$C$ is special in the sense of Griffiths-Harris \cite {GH} page
251, because $K_C - D = (1 + (n-2) r) K_X | C \ge (r+1) K_X | C
\ge 0$. Applying Clifford's Theorem [ibid], one has $h^0(C, D) \le
1 + \deg(D)/2$. Substituting this and $\deg(D) = r K_X .
(rK_X)^{n-1} = r^n K_X^n$ into the above, we obtain (4).
\end{proof}

In the case of generically finite surjective map between varieties
of general type, the volume upstairs can control the volume
downstairs and sometimes even the degree of the map. See also
Lemma \ref{Kdeg'}.

\begin{lem}\label{Kdeg}
Let $X$ and $Y$ be smooth projective $n$-folds of general type and
$f : X \rightarrow Y$ a generically finite surjective morphism.

\begin{list}{}{
\setlength{\leftmargin}{10pt}
\setlength{\labelwidth}{6pt}
}
\item[\rm{(1)}]
We have $V(X) \ge \deg(f) \,\, V(Y)$.
\item[\rm{(2)}]
Suppose further that $|K_Y|$ defines a generically finite rational
map $\Phi$. Then $V(Y) \ge \deg \Phi \,\, \deg \Phi(Y) \ge 1$ and
$\deg(f) \le V(X)$.
\end{list}
\end{lem}

\begin{proof}
Let $\varepsilon > 0$. By Fujita's approximation as in Lazarsfeld
\cite{La} Theorem 11.4.4, after smooth modifications of $X$ and
$Y$, we may assume that $K_Y \sim_{\bf Q} H + E$ with $H$ an ample
${\bold Q}$-divisor and $E$ an effective ${\bold Q}$-divisor, such
that $H^n = V(H) > V(Y) - \varepsilon$. By the ramification
divisor formula, we have $V(X) \ge V(f^*H) = (f^*H)^n = \deg(f)
H^n > \deg(f) (V(Y) - \varepsilon)$. Now (1) follows by letting
$\varepsilon$ tend to zero.

The first inequality in (2) is by the proof of Hacon-McKernan
\cite{HM-birat} Lemma 2.2; the rest now follows.
\end{proof}

The assertion (1) below means that the volume controls the order of
a free-acting group, while (2) says that sometimes $|G|$ can be
bounded by its exponent $\exp(G)$. We denote by $\cartier(Y)$ the
Cartier index, i.e., the smallest positive integer such that $c K_Y$
is a Cartier divisor.

\begin{lem}\label{Kquot}
Let $X$ be a smooth minimal projective $n$-fold of general type
and $f : X \rightarrow Y$ a finite surjective morphism.

\begin{list}{}{
\setlength{\leftmargin}{10pt}
\setlength{\labelwidth}{6pt}
}
\item[\rm{(1)}]
Suppose that $f$ is unramified. Then $\deg(f) = K_X^n/K_Y^n \le
K_X^n$. In particular, if a subgroup $H \le \Aut(X)$ acts freely
on $X$, then $|H| = K_X^n / K_Y^n \le K_X^n$, where $Y = X/H$.

\item[\rm{(2)}]
Suppose that a subgroup $G \le \Aut(X)$ has isolated fixed locus
$X^G$. Let $f : X \rightarrow Y := X/G$ be the quotient map and $c
= \cartier(Y)$ the Cartier index. Then $c \, | \, \exp(G)$ and
$|G| \le c K_X^n \le \exp(G) \,\, K_X^n$.
\end{list}
\end{lem}

\begin{proof}
We note that $K_X$ is nef and big and that $K_X^n$
is a positive integer.

(1) follows from the fact that $K_X = f^*K_Y$.

(2) Denote by $G_x = \{g \in G \, | \, g(x) = x\}$ the stabilizer
subgroup at the point $x \in X$. We may regard $G_x$ as a subgroup
of $\GL_n(T_{X, x}) \cong \GL_n({\bold C})$. For the image
$\overline{x} \in Y$ of $x$, the local Cartier index $c({\overline
x})$ equals the index of the quotient group $G_x/G_x \cap
\SL_n(\CC) \subset \GL_n(\CC)/\SL_n(\CC) \cong \CC^*$, which is
cyclic and generated by the image $\overline{g}$ of some $g \in
G_x$. Thus $c({\overline x})$ divides $\exp(G_x)$ and also
$\exp(G)$. Hence $c = \cartier(Y)$, which is the $\lcm$ of
$c({\overline x})$, also divides $\exp(G)$. Since $c K_Y$ is
Cartier, it can be written as $H_1 - H_2$ with very ample divisors
$H_j$ not passing through the isolated set $\Sing(Y)$ (the image
of $X^G$). Then $c K_Y^n = (H_1 - H_2) . K_Y^{n-1} = (K_{Y} |
H_1)^{n-1} -  (K_Y | H_2)^{n-1}$ is an integer. Since $f : X
\rightarrow Y$ is etale outside the finite set $\Sing(Y)$, one has
$K_X = f^*K_Y$. Hence $|G| = \deg(f) = K_X^n / K_Y^n = c K_X^n / c
K_Y^n \le c K_X^n$.
\end{proof}

An application of the result below is given in Remark \ref{remBb}.
We remark that the finiteness assumption below on $X^g$ can be
removed for surfaces as shown in Ueno \cite{Ue2}.

\begin{lem} \label{Lefix}
Suppose that $X$ is a smooth projective $n$-fold and $g$ is an
automorphism of finite order. Assume that $X^g$ is finite.

\begin{list}{}{
\setlength{\leftmargin}{10pt}
\setlength{\labelwidth}{6pt}
}
\item[\rm{(1)}]
The Euler number of the fixed locus $X^g$ is given by $e(X^g) =$
$\sum_{i=0}^{2n}$ $(-1)^i \Tr \ g^* | H^i(X, \ZZ)/(\torsion)$.
\item[\rm{(2)}]
One has $|e(X^g)| \le \sum_{i=0}^{2n} B_i(X)$, where $B_i(X)$ is the
Betti number.

\end{list}
\end{lem}

\begin{proof}
(1) is the classical Lefschetz fixed point formula. For (2), let
$\ord(g) = m$ and diagonalize $g^* | H^i(X, \CC) =
\diag[\zeta_m^{s_1}, \dots, \zeta_m^{s_{b_i}}]$ where $b_i = B_i(X)$
and $\zeta_m^{s_j}$ is an $m$-th root of $1$. Thus $|\Tr g^* |
H^i(X, \ZZ)/\torsion| = |\sum_{j=1}^{b_i} \zeta_m^{s_j}| \le b_i$.
This proves the lemma.
\end{proof}

The following is the relation between two stabilizers
$G_x$ and $(G/H)_{\overline x}$.

\begin{lem}\label{GGx}
Suppose that a normal subgroup $H \trianglelefteq G$ acts
freely on a variety $X$. Set $\overline{G} = G/H$
and $\overline{X} = X/H$.

\begin{list}{}{
\setlength{\leftmargin}{10pt}
\setlength{\labelwidth}{6pt}
}
\item[\rm{(1)}]
For every $x \in X$ and its image $\overline{x} \in \overline{X}$,
the homomorphism between stabilizers $\varphi : G_x \rightarrow
G_{\overline x}$  which is given by $g \mapsto \overline{g} = gH$,
is an isomorphism.
\item[\rm{(2)}]
In particular, the fixed locus $X^G$ is empty (resp. finite) if
and only if $\overline{X}^{\overline{G}}$ is empty (resp. finite).
\end{list}
\end{lem}

\begin{proof}
To show that $\varphi$ is surjective, suppose that $\overline{g}
(\overline{x}) = \overline{x}$. Then $g(x) = h(x)$ for some $h \in
H$. Thus $h^{-1} g \in G_x$ and $\overline{g} = \varphi(h^{-1}g)$.
To show that $\Ker(\varphi)$ is trivial, suppose that $\varphi(g)
= \id \in \overline{G}_{\overline x} \le \overline{G}$. So $g \in
H$. Hence $g \in H \cap G_x = \{\id\}$ by the freeness of the $H$
action on $X$. So $\Ker(\varphi)$ is trivial and the lemma is
proved.
\end{proof}

Next are about nilpotent groups and relation between centralizers.

\begin{lem} \label{gTh}
Let $G$ be a finite group.

\begin{list}{}{
\setlength{\leftmargin}{10pt}
\setlength{\labelwidth}{6pt}
}
\item[\rm{(1)}]
If $G$ is nilpotent, then $\exp(G)$ equals $\mu(G) := \max$
$\{\ord(g) \, | \, g \in G\}$.
\item[\rm{(2)}]
Consider quotient groups inclusion $K/H \trianglelefteq G/H$. Let
$\tau \in K$ and $\overline{\tau} = \tau H \in G/H$. Then
$|C_{K/H}({\overline \tau})| \le |C_K(\tau)| \le |C_G(\tau)|$.
\end{list}
\end{lem}

\begin{proof}
(1) is true for $p$-groups while $G$, being nilpotent, is a direct
product of its Sylow subgroups (see Gorenstein \cite{Go} Theorem
3.5). So (1) follows.

(2) Write $C_{K/H}({\overline \tau}) = L/H$. Consider the map:
$\varphi : L \rightarrow \tau H$ ($g \mapsto g^{-1} \tau g$). Then
$|L/C_L(\tau)| = |\Imm(\varphi)| \le |H|$, so $|L/H| \le
|C_L(\tau)|$.
\end{proof}

When $\Aut(X)$ fixes a point $x$, Jordan's lemma and its
consequences below, together with Xiao's Theorem \ref{XiaoTh},
will produce a linear bound $|\Aut(X)| \le J_nx_n K_X^n$.

\begin{lem} \label{Jordan} The following are true.
\begin{list}{}{
\setlength{\leftmargin}{10pt}
\setlength{\labelwidth}{6pt}
}
\item[\rm{(1)}] ({\bf Jordan's lemma})
Let $J_n = (n+2)!$ (resp. $J_n = n^4 (n+2)!$) for $n > 63$ (resp.
for $n \le 63$). Then every finite group in $\GL_n(\CC)$ contains
an abelian normal subgroup of index $\le J_n$.
\item[\rm{(2)}]
Suppose that $X$ is an $n$-fold and $G \le \Aut(X)$ is a finite subgroup
fixing a smooth point $x \in X$. Then $G$ contains an abelian normal
subgroup of index $\le J_n$.
\item[\rm{(3)}]
Suppose that $X$ is a projective $n$-fold of general type and $x
\in X$ a smooth point. If $G \le \Aut(X)$ fixes $x$, then $G$
contains an abelian normal subgroup of index $\le J_n$.
\end{list}
\end{lem}

\begin{proof}
For (1), see Weisfeiler \cite{We} page 5279 (for better bound) and
Aljadeff- Sonn \cite{AS} page 353. The original number produced by
Jordan was $J(n) = (49n)^{n^2}$. For (2), $G$ can be regarded as a
subgroup of $\GL_n(T_{X, x}) \cong \GL_n(\CC)$, so apply (1). For
(3), $\Aut(X)$ and hence $G$ are finite since $X$ is of general
type. Then apply (2).
\end{proof}

The result below has been proved by many authors. We only mention
that the nefness of $\Omega_X^1$ comes from the exact sequence
below and the fact that quotient bundles of a nef bundle are again
nef (see Hartshorne \cite{Hart} Ch II, Theorem 8.17; Lazarsfeld
\cite{La} Proposition 6.1.2):
$$0 \rightarrow N^{\vee}_{X/A} \rightarrow
\Omega_A^1 | X \rightarrow \Omega_X^1
\rightarrow 0.$$

\begin{thm} \label{Gauss}
Let $A$ be an abelian variety and $X \subset A$
a subvariety of general type.

\begin{list}{}{
\setlength{\leftmargin}{10pt}
\setlength{\labelwidth}{6pt}
}
\item[\rm{(1)}]
$|K_{\tilde X}|$ defines a generically finite map, where
$\tilde{X} \rightarrow X$ is any desingularization.

\item[\rm{(2)}]
Suppose that $X$ is smooth. Then
$\Omega_X^1$ is nef and $K_X$ is ample.
Moreover, $|K_X|$ is base point free
and $\Phi_{|K_X|}$ is a finite morphism.
\end{list}
\end{thm}

\begin{proof}
For the proof, see Ran \cite{Ra} Corollary 2, or Abramovich
\cite{Ab}. See also Hartshorne \cite{Hart71}, Ueno \cite{Ue} and
Griffiths-Harris \cite{GH79} for earlier development.
\end{proof}

Here are some properties of simple abelian varieties to be used in
\S 6.

\begin{lem}\label{sA}
Let $A$ be a simple abelian variety.

\begin{list}{}{
\setlength{\leftmargin}{10pt}
\setlength{\labelwidth}{6pt}
}
\item[\rm{(1)}]
Every proper subvariety of $A$ is of general type.
\item[\rm{(2)}]
If $\id \ne g \in \Aut_{\variety}(A)$ is not a translation, then
$A^g$ is a non-empty finite set.
\item[\rm{(3)}]
If $X \rightarrow A$ is a non-contant morphism from a projective
variety $X$, then $q(\tilde{X}) \ge \dim(A)$, where $\tilde{X}
\rightarrow X$ is any desingularization.
\end{list}
\end{lem}

\begin{proof}
(1) is just Ueno \cite{Ue} Corollary 10.10.

(2) Suppose that $g$ is not a translation. By Birkenhake-Lange
\cite{BL} Lemma 13.1.1, we may assume that $\id \ne g \in
\Aut_{\group}(A)$.  Now $\Ker(g - \id)$ is a proper subgroup of $A$,
whence it is non-empty and $0$-dimensional by the simplicity of $A$.
So (2) follows (see \cite{BL} Formula 13.1.2).

(3) By the universal property of $\Alb(\tilde{X})$, the
non-constant morphism $\tilde{X} \rightarrow A$ factors as
$\tilde{X} \rightarrow \Alb(\tilde{X}) \rightarrow A$, where the
latter map is a homomorphism (modulo a translation) and has image
the translation of a non-trivial subtorus of $A$. This image is
$A$ by the simplicity of $A$. This implies (3) :  $\dim(A) \le
\dim \Alb(\tilde{X}) = q(\tilde{X})$.
\end{proof}

By verifying the case $\ell = 8$ (resp. $\ell = 12$) and
using induction, we have (1) (resp. (2)) below to be used in Section 5.

\begin{lem}\label{alt}
Let $A_k$ be the alternating group of order $k!/2$.
Let $\ell = [k/4] = \max\{z \in \ZZ \, | \, z \le k/4\}$
so that $4 \ell \le k \le 4 \ell + 3$.

\begin{list}{}{
\setlength{\leftmargin}{10pt}
\setlength{\labelwidth}{6pt}
}
\item[\rm{(1)}]
If $\ell \ge 8$, then $|A_k| \le (4 \ell +3)!/2 \le |A_{3 \ell}|^{1.7}$.
\item[\rm{(2)}]
If $\ell \ge 12$, then $|A_k| \le |A_{\ell}|^8$.
\end{list}
\end{lem}

We end \S 2 with Xiao's linear bound for abelian subgroups:

\begin{thm} \label{XiaoTh} (Xiao \cite{Xi96} Theorem 1)
Let $X$ be a smooth projective $n$-fold with nef and big $K_X$. Let
$G \le \Aut(X)$ be an abelian subgroup. Then there exists a constant
$x_n$, independent of $X$, such that $|G| \le x_n K_X^n$.
\end{thm}
\section{\bf G-equivariant fibrations; the proof of Theorem \ref{Thconj}}
In this section we consider $G$-equivariant fibration $f : X
\rightarrow Y$. We will prove simultaneously Theorem \ref{Th3.1}
below and Theorem \ref{Thconj} in the Introduction.

\begin{thm}\label{Th3.1}
Let $X$ and $Y$ be smooth projective varieties of general type and
dimensions $n$ and $k$, respectively. Suppose that a group $G$ acts
faithfully on $X$ and acts anyhow on $Y$. Let $f : X \rightarrow Y$
be a $G$-equivariant surjective morphism with connected general
fibre $F$. Suppose further that $|\Bir(Y)| \le a_1 \ V(Y)^{d}$ and
$|\Bir(F)| \le a_2 \ V(F)^{d}$ for some positive integers $a_i$ and
$d$.
\par
Then we have:
$$|G| \le a \, V(X)^d, \hskip 2pc \text{\rm with} \hskip 1pc a = \frac{a_1
a_2}{(\binom{n}{k})^d} \,\, .$$
\end{thm}

\par
We now prove first Theorem \ref{Th3.1} and late Theorem \ref{Thconj}
as an application. Note that the general fibre $F$ is also of
general type since the Kodaira dimension $\kappa(X) = \dim X$ and by
Iitaka's easy addition: $\kappa(X) \le \kappa(F) + \dim Y$.
 Since $f$ is $G$-equivariant, we have the
natural exact sequence:
$$1 \rightarrow K \rightarrow G \rightarrow \Aut(Y).$$
The group $K$ acts trivially on $Y$ and hence faithfully on the
fibre $F$, so $K$ can be regarded as a subgroup of $\Aut(F)$, where
$F$ is smooth and of general type. Now Theorem \ref{Th3.1} follows
from the Theorem in the Appendix:
$$
|G| \le |K| \, |\Aut(Y)| \le |\Aut(Y)| \ |\Aut(F)| \le$$
$$|\Bir(Y)| \ |\Bir(F)|
\le a_1 a_2 \ V(Y)^d \ V(F)^d \le a V(X)^d.
$$
\par
Finally, Theorem \ref{Thconj} is a consequence of Theorem
\ref{Th3.1} and the inductive assumption, where we set $a_1 =
\delta_{\dim Y}$, $a_2 = \delta_{\dim F}$ and $d = 1$.

\begin{rem}
Weaker versions of Theorems \ref{Thconj} and \ref{Th3.1} can also
be proved using Lemma \ref{smalldiminv}, Kollar \cite{KoI} Theorem
3.5 (ii) and Catanese-Schneider \cite{CS} pages 10-11, which are
strengthened to the current form, thanks to the Appendix - an
application of the weak positivity due to Fujita \cite{Fuj},
Kawamata \cite{Ka82b} and Viehweg \cite{Vi82b}.
\end{rem}
\section{\bf Bound of Betti numbers; the proof of Theorem \ref{Bbetti}}
In this section, for a smooth projective variety $X$, we shall bound
the invariants like Betti numbers $B_i(X)$ and the Euler number
$e(X)$, in terms of the volume of $X$. We first bound the
irregularity and geometric genus.

\begin{prop} \label{qp}
Let $X$ be a smooth projective $n$-fold ($n \ge 2$) with
ample $K_X$. Then $p_g(X) \le e_n K_X^n$
and $q(X) \le e_n K_X^n$, where
$e_n = \max\{n + r^n, \, 2 + \frac{1}{2} r^{n-2} (1 + (n-2) r)^2\}$
with $r = 2 + n(n+1)/2$.
\end{prop}

\begin{proof}
For $r = 2 + n(n+1)/2$, as in Lemma \ref{smalldiminv}, $|rK_X|$ is
base point free. Since $K_X$ is ample, $\Phi = \Phi_{|rK_X|} : X
\rightarrow \PP^N$ with $N = P_r(X) - 1$, is a finite morphism.
Hence $(rK_X)^n = \deg(\Phi) \deg \Phi(X) \ge \deg \Phi(X) \ge N -
n + 1$ (see Griffiths-Harris \cite{GH} page 173). Thus $p_g(X) \le
P_r(X) = N + 1 \le n + r^n K_X^n \le (n + r^n) K_X^n$.

As in Lemma \ref{smalldiminv}, take the ladder
$X = X_{n} \supset X_{n-1} \supset \cdots \supset X_2$,
where $X_i$ is the intersection of $n-i$ general members of $|rK_X|$,
and a smooth $i$-fold with ample $K_{X_i} = (1 + (n-i)r) K_X | X_i$.
Each $X_{i-1}$ is a smooth ample divisor on $X_i$.
By Lefschetz hyperplane section theorem,
we have $q(X) = q(X_{n-1}) = \dots = q(X_2)$.
By Lemma \ref{smalldiminv}, we conclude the proposition:
$q(X_2) \le \frac{1}{2}K_{X_2}^2 + 2
= 2 + \frac{1}{2} r^{n-2} (1 + (n-2) r)^2 K_X^n$.
\end{proof}

To bound the Betti numbers, we need the bound of Euler number $e(X)$ below.

\begin{lem}\label{eX}
Let $X$ be a smooth projective $n$-fold.

\begin{list}{}{
\setlength{\leftmargin}{10pt}
\setlength{\labelwidth}{6pt}
}
\item [\rm{(1)}]
If the cotangent bundle $\Omega_X^1$ is nef,
then $0 \le (-1)^n e(X) = c_n(\Omega_X^1) \le K_X^n$.
\item [\rm{(2)}]
If $K_X$ is ample then $|e(X)| \le h_n K_X^n$, where $h_n$ is the
constant, independent of $X$, in Heier \cite{He} Proposition 1.7.
\end{list}
\end{lem}

\begin{proof}
(2) is proved in Heier \cite{He} Proposition 1.7. For (1), by the
nefness of $\Omega_X^1$ and Demailly-Peternell-Schneider \cite{DPS}
page 31, we have $0 \le c_n(\Omega_X^1) \le c_1(\Omega_X^1)^n =
K_X^n$. Note also that $c_n(\Omega_X^1) = (-1)^n c_n(T_X) = (-1)^n
c_n(X) = (-1)^n e(X)$, where $T_X$ is the tangent bundle; see for
instance, Fulton \cite{Ful} Example 3.2.13. This proves the lemma.
\end{proof}

\begin{setup}
{\bf Proof of Theorem \ref{Bbetti}}
\end{setup}

We will prove by induction on $n$. The case $n = 1$ is clear. When
$n = 2$, it follows from Lemma \ref{smalldiminv} or \ref{qp}.
Indeed, note that $B_1(X) = 2q(X)$, $\chi(\OO_X) \le 1 + p_g(X)$,
$c_2(X) \le K_X^2 + c_2(X) = 12 \chi(\OO_X)$ and $B_2(X) = c_2(X)
- 2 + 4q(X)$.

Suppose $n \ge 3$ and the theorem is true for all such $k$-folds
with $2 \le k < n$. Let $X_{n-1} $ be a general member of $|rK_X|$
with $r = 2 + n(n+1)/2$ (see Lemma \ref{smalldiminv}). Note that
$K_{X_{n-1}} = (1 + r)K_X | X_{n-1}$ is ample and
$K_{X_{n-1}}^{n-1} = r(1 + r)^{n-1} K_X^n$. By induction,
$B_i(X_{n-1}) \le a_{n-1} K_{X_{n-1}}^{n-1} = a_n' K_X^n$ where
$a_n' = a_{n-1} r (1+r)^{n-1}$. By Lefschetz hyperplane section
theorem, $B_i(X) = B_i(X_{n-1})$ for all $i \le n-2$ and
$B_{n-1}(X) \le B_{n-1}(X_{n-1})$ (see Lazarsfeld \cite{La}
Theorem 3.1.17). Note also that $B_j(X) = B_{2n-j}(X)$ by
Lefschetz duality. Thus $B_{j}(X) \le B_{j}(X_{n-1}) \le a_n'
K_X^n$ for all $j \ne n$.

On the other hand, by Lemma \ref{eX}, $h_n K_X^n \ge
|\sum_{i=0}^{2n} (-1)^{i} B_i(X)|$.
So $|B_n(X)| \le h_n K_X^n + |\sum_{i \ne n} (-1)^{i} B_i(X)|
\le a_n K_X^n$ where $a_n = h_n + 2n a_n'$.
We are done by the induction.
This proves Theorem \ref{Bbetti}.

\begin{rem} \label{remBb}
For the $X$ in Theorem \ref{Bbetti}, it might be possible to give a
second proof of Xiao's linear bound (in terms of $K_X^n$) of
$\ord(g)$ for every $g \in \Aut(X)$ with $\ord(g)$ prime and $X^g$
finite. Indeed, note that the quotient map $X' = X - X^g \rightarrow
X'/\langle g \rangle$ is unramified, and hence $\ord(g)$ divides the
Euler number $e = e(X) - e(X^{g})$, so $e$ is bounded by $|e(X)| +
|e(X^g)|$ (provided that $e \ne 0$) while the latter is bounded
linearly by $K_X^n$ (see Lemmas \ref{Lefix} and \ref{eX} and Theorem
\ref{Bbetti}).
\end{rem}
\section {\bf Automorphism groups of subvarieties of abelian varieties}
In this section, we shall prove the following two results:

\begin{thm} \label{Th2'}
Let $X$ be a smooth $n$-fold ($n \ge 3$) of general type contained
in an abelian variety $A$ of dimension $q$. Let $G$ be a subgroup of
$\{g \in \Aut_{\variety}(A) \, | \, g(X)$ $= X\}$ such that the
fixed locus $A^g$ for every $\id \ne g \in G$ is a non-empty finite
set unless $g$ is a translation of $A$. Set $V = K_X^n$ (see Theorem
\ref{Gauss}).

Then there is a constant $d_n$, independent of $X$ and $A$, such
that $|G| \le d_n q^{10-b}V^b \le  d_n S^{10}$, where $S = \max\{V,
q\}$ and $5 \le b \le 10$.
\end{thm}

\begin{cor}\label{Th2''}
Let $A$ be a simple abelian variety and
$X \subset A$ a smooth $n$-dimensional ($n \ge 3$) proper subvariety.
Let $G = \{g \in \Aut_{\variety}(A) \, | \,$ $g(X)= X\}$ be the stabilizer.
Set $V = K_X^n$ (see Theorem \ref{Gauss}).
Then there is a constant $d_n'$ (independent of $X$ and $A$) such that
$|G| \le d_n' V^{10}$.
\end{cor}

We fix an abelian variety $A$ of dimension $q$ and an
$n$-dimensional smooth subvariety $X \subset A$ of general type. Let
$G$ be a subgroup of the group $\{g \in \Aut_{\variety}(A) \, | \,
g(X) = X\}$. Note that $\Aut_{\variety}(A) = T  \rtimes A_0$ (split
extension) where $T = \{T_t \, | \, t \in A\}$ is the normal
subgroup of translations and $A_0 = \Aut_{\group}(A)$ is the
subgroup of bijective group-homomorphisms.

\begin{setup}\label{setup1'}
{\bf Assumption}.
\begin{list}{}{
\setlength{\leftmargin}{10pt}
\setlength{\labelwidth}{6pt}
}
\item [\rm{(i)}]
$X \subset A$ is a smooth projective $n$-fold of general type in the
abelian variety $A$ of dimension $q$,

\item [\rm{(ii)}]
the subgroup $G \subseteq \{g \in \Aut_{\variety}(A) \, | \, g(X) = X\}$
contains no translations of $A$
(we note that $G \le \Aut(X)$ is finite
because $X$ is of general type), and

\item [\rm{(iii)}]
the fixed locus $A^g$ for every $\id \ne g \in G$ is a non-empty
finite set.
\end{list}
\end{setup}

We collect some information on the structure of $G$.

\begin{lem}\label{gr}
Suppose that $G \subseteq \Aut(X)$ and $X \subset A$ satisfy
conditions in \ref{setup1'}. Let $\id \ne g \in G$. Then we have:

\begin{list}{}{
\setlength{\leftmargin}{10pt}
\setlength{\labelwidth}{6pt}
}
\item [\rm{(1)}]
Let $G_0$ be the image of the composition $G \subseteq \Aut_{\variety}(A)
\rightarrow A_0$. Then $G \rightarrow G_0$ is an isomorphism.
\item [\rm{(2)}]
For $g \in G$ and $g_0$ its image in $G_0$, we have
$|A^g| = |A^{g_0}| \ge |\{0\}| = 1$.
\item [\rm{(3)}]
If $H \le G$ is abelian, then $H$ is cyclic.
\item [\rm{(4)}]
For $H \le G$, if $|H| = p^2$ for some prime $p$, then $H$ is cyclic.
\item [\rm{(5)}]
If $|X^g| \ge 2$ for some $g \in G$, then $\ord(g) = p^s$ for some prime $p$.
\item [\rm{(6)}]
If $g$ has order $m$ in $G$ and if $\zeta_m$ denotes a primitive
$m$-th root of $1$, then there is is a diagonalization $g^* | H^0(A,
\Omega_A^1)$ $= \diag[\zeta_m^{s_1}, \dots, \zeta_m^{s_q}]$ where
each $\zeta_m^{s_j}$ is a primitive $m$-th root of $1$.
\item [\rm{(7)}]
If $g \in G$, then the Euler function $\varphi(\ord(g))$ divides
$2q$.
\item [\rm{(8)}]
If $p^2$ divides $|G|$ for some prime $p$, then there is an
element $g \in G$ such that $\ord(g) = p^2$.
\item [\rm{(9)}]
If $g \in G$, then $\ord(g)$ divides $4q^2 \alpha(g)$, where
$$\alpha(g) = \prod_{\text{\rm prime}
\, p \,|\, \ord(g), \, p^2 \nmid |G|} \, p.$$
\item [\rm{(10)}]
If $H \le G$, then the exponent $\exp(H)$ divides $4q^2
\alpha(H)$, where
$$\alpha(H) = \prod_{\text{\rm prime}
\, p \,|\, |H|, \, p^2 \nmid |G|} \, p.$$
\end{list}
\end{lem}

\begin{proof}
\ref{setup1'} (ii) implies (1). Here (2), (3), (6) and (7) are
consequence of \ref{setup1'} (iii) and Birkenhake - Lange \cite{BL}
Lemma 13.1.1 and Propositions 13.2.2, 13.2.1 and 13.2.5. Our (4) is
because a group of order $p^2$ is abelian and by (3). Now (5)
follows from (1) and Birkenhake - Lange \cite{BL} Corollary 13.2.4.
Also (8) is a consequence of (4) together with Sylow's theorem and a
basic result on $p$-groups. For (9), write $\ord(g) = p_1^{t_1}
p_2^{t_2} \dots p_u^{t_u}$ with $t_i \ge 1$. Then $\varphi(\ord(g))
= p_1^{t_1-1}(p_1 - 1) \dots p_u^{t_u - 1}(p_u - 1)$. If $t_i \ge
2$, then $2(t_i - 1) \ge t_i$, and $p_i^{t_i}$ divides
$(\varphi(\ord(g)))^2$ as well as $(2q)^2$ by (7). If $p^2 \, | \,
|G|$, then $p \, | \, 2q$ by (8) and (7). Thus $\ord(g) \, | \, 4q^2
\alpha(g)$ and (9) is proved. Our (9) implies (10). This proves the
lemma.
\end{proof}

The centralizer lemma below is a key to the proof of Theorem \ref{Th2'}.

\begin{lem}\label{central} ({\bf centralizer lemma})
Suppose that $G \subseteq \Aut(X)$ and $X \subset A$ satisfy
conditions in \ref{setup1'}. Let $\tau \in G$ be of prime order $p$
such that the fixed locus $X^{\tau}$ is non-empty. Let $V = K_X^n$
(see Theorem \ref{Gauss}).
\begin{list}{}{
\setlength{\leftmargin}{10pt}
\setlength{\labelwidth}{6pt}
}
\item [\rm{(1)}]
Every prime factor $p_1 \ne p$ of the order $|C_G(\tau)|$ of the
centralizer, divides $|X^{\tau}| (|X^{\tau}| - 1)$.
\item [\rm{(2)}]
$\exp(C_G(\tau))$ divides $4 p^{\varepsilon} q^2|X^{\tau}|
(|X^{\tau}| - 1)$, where $\varepsilon = 0$ (resp. $1$) when $p^2$
divides $|G|$ (resp. otherwise).
\item [\rm{(3)}]
One has $|C_G(\tau)| \le k_n q^2V^{3 + \varepsilon}$, where $k_n =
\max\{J_n x_n, \, 4a_n^2 x_n^{\varepsilon}\}$; see Lemma
\ref{Jordan} and Theorems \ref{Bbetti} and \ref{XiaoTh}.
\end{list}
\end{lem}

\begin{proof}
(1) Suppose $g \in C_G(\tau)$ has order equal to a prime $p_1 \ne p$.
Write $|X^{\tau}| \equiv r$ with $0 \le r < p_1$. If $r = 0, 1$,
then we are done.
Suppose that $r \ge 2$. Then $|X^{\tau g}| \ge r \ge 2$.
By Lemma \ref{gr}, $p p_1 = \ord(\tau g)$ equals a prime power,
a contradiction.

(2) follows from (1) and Lemma \ref{gr} (10).

(3) If $|X^{\tau}| = 1$, then $C_G(\tau)$ fixes the unique point
$x$ in $X^{\tau}$. Hence $C_G(\tau) \le \GL_n(T_{X, x})$ and we
apply Lemma \ref{Jordan} and Theorem \ref{XiaoTh} and obtain
$|C_G(\tau)| \le J_n x_n K_X^n$. Suppose now that $|X^{\tau}| \ge
2$. By (2), Lemmas \ref{Kquot} and \ref{Lefix} and Theorem
\ref{Bbetti}, we have $|C_G(\tau)| \le \exp(C_G(\tau) V \le 4
p^{\varepsilon} q^2|X^{\tau}| (|X^{\tau}| - 1) V \le 4
p^{\varepsilon} q^2 (\sum_{i=0}^{2n} B_i(X))^2 V \le 4a_n^2
p^{\varepsilon} q^2 V^3$. By Theorem \ref{XiaoTh}, $p \le x_n V$
and hence (3) follows. This proves the lemma.
\end{proof}

Theorem \ref{Th2'} should follow from the crucial proposition below.

\begin{prop} \label{key}
Suppose that $X \subset A$ and $G \subseteq \Aut(X)$ satisfy the
conditions in \ref{setup1'}. Let $q = \dim(A)$ and $V = K_X^n$ (see
Theorem \ref{Gauss}). Then the there is a constant $d_n$,
independent of $X$ and $A$, such that $|G| \le d_n q^{10-b} V^b \le
d_n S^{10}$, where $S = \max\{V, q\}$ and $5 \le b \le 10$.
\end{prop}

For the proof, we argue closely along the lines of
Huckleberry-Sauer \cite{HS} and Szabo \cite{Sz} (see \cite{As}).
But we use the centralizer Lemma \ref{central} instead. Let $H
\trianglelefteq G$ be a maximum among normal subgroups of $G$
which acts freely on $X$. Take a normal subgroup $K/H
\trianglelefteq G/H$. Fix some $\tau \in K$ such that the fixed
locus  $X^{\tau} \ne \emptyset$. We may assume that $\ord(\tau) =
\ord(\overline{\tau}) = p$ is a prime (see \ref{setup1'} (iii) and
Lemma \ref{GGx}).

\begin{claim} \label{GG1}
Let $G_1 = C_G(\tau)$. Then $|G| \le |G_1| |K|$.
\end{claim}

\begin{proof}
Note that $\varphi : G \rightarrow K$ ($g \mapsto g^{-1} \tau g$)
is a well-defined map.
Then the claim follows from that $|G/G_1| = |\varphi(G)| \le |K|$.
\end{proof}

Suppose that we can choose $K/H$ to be abelian.
Note that $K/H$ acts on the smooth minimal $n$-fold $X/H$ of general type.
By Theorem \ref{XiaoTh}, we have $|K/H| \le x_n K_{X/H}^n
= x_n V/|H|$. Thus by Claim \ref{GG1} and
Lemmas \ref{central} and \ref{Kquot},
we have $|G| \le |G_1| |K/H| |H| \le k_nq^2 V^{3 + \varepsilon} x_n V
\le d_n q^2 V^{5}$ with $d_n = k_n x_n$.

We may assume that there is no such abelian $K/H$.
Then $G/H$ is semi-simple, i.e., it has no non-trivial
abelian normal subgroup.
Let $M/H$ be a non-trivial minimal normal subgroup. Then
$M/H = \prod^k E/H$, a $k$-fold product of the same
non-abelian simple group $E/H$.
Let $M_i/H$ be all distinct non-trivial minimal normal subgroups of $G/H$
and let $S/H = \prod M_i/H$ be the sockel of $G/H$.
We write $M_i/H = \prod^{k_i} E_i/H$ with non-abelian simple group $E_i/H$.

Suppose that $M/H = \prod^k E/H$ with $k \ge 2$. Then $\alpha(M)
\, | \, |H|$ in notation of Lemma \ref{gr}, because for every
prime factor $p_1$ of $|M|$, either $p_1 \, | \, |H|$, or $p_1 \,
| \, |M/H|$  and hence $p_1^k \, | \, |M/H|$ (and $p_1^2 \, | \,
|G|$). So by Lemma \ref{gr}, $\exp(M)$ divides $4q^2 \alpha(M)$
and also $4q^2 |H|$. Our Lemma \ref{Kquot} implies that $|M| \le
\exp(M) K_X^n$. Substituting all these in and applying Lemmas
\ref{central} and \ref{Kquot}, we obtain $|G| \le |G_1| |M| \le
k_nq^2 V^{3 + \varepsilon} 4q^2 |H| V \le 4k_n q^4 V^6$ and we are
done.

So we may assume that $M/H = E/H$ is non-abelian simple
for every non-trivial minimal normal $M/H \trianglelefteq G/H$.
Suppose that $M/H$ is one of the $26$ sporadic non-abelian simple groups.
Then $M/H \le d$, a constant (the order of the Monster simple group).
Thus as above,
$|G| \le |G_1| |M/H| |H| \le k_nq^2 V^{3 + \varepsilon}
d V \le d k_n q^2 V^5$.

Suppose that $E/H$ is of Lie type. By the proof of
Huckleberry-Sauer \cite{HS} Proposition 7, there exist a universal
constant $d$ and a Sylow $p_1$-subgroup $U/H$ of $M/H$ such that
$|M/H| \le d |U/H|^{5/2}$. Note that $U/H$ acts on the smooth
minimal $n$-fold $X/H$ of general type and the fixed locus
$(X/H)^{U/H}$ is finite by \ref{setup1'} (iii) and Lemma
\ref{GGx}. By Lemmas \ref{Kquot} and \ref{gTh} and Theorem
\ref{XiaoTh}, we have $|U/H| \le \exp(U/H) K_{X/H}^n = \mu(U/H)
V/|H| \le x_n K_{X/H}^n V/|H| = x_n (V/|H|)^2$. Thus $|M/H| \le d
x_n^{5/2}$ $(V/|H|)^5$. By Lemma \ref{central} we have $|G| \le
|G_1| |M| \le k_nq^2 V^{3 + \varepsilon} d x_n^{5/2} V^5$ $= d_n
q^2 V^{8 + \varepsilon}$ with $d_n = k_n d x_n^{5/2}$. If $p^2 \,
| \, |G|$ then $\varepsilon = 0$ and we are done.

If the fixed locus $X^U = \emptyset$, then $|U| \le V$ by Lemma \ref{Kquot}
and we will even have a better bound.
We may assume that the fixed locus $X^U \ne \emptyset$
(equivalently $(X/H)^{U/H} \ne \emptyset$ by Lemma \ref{GGx}).
Then our initial $\overline{\tau}$ (of order $p$) can be taken from $U/H$
so that $p_1 = p$ and $U/H$ is a $p$-Sylow subgroup of $M/H$.
If $p^2$ divides $|U/H|$ then it also divides $|G|$,
whence $\varepsilon = 0$
and we are done. Otherwise, $|U/H| = p \le x_n K_{X/H}^n = x_n V/|H|$
by Theorem \ref{XiaoTh} and we will be done again.

We are left with the case where each $M_i/H$ is an alternating
group $A_{k_i}$. Suppose that $M/H = A_k$ is the smallest among
them. If there are two $M_i/H$ then we will be done as in the case
$M/H = \prod^k E/H$ with $k \ge 2$ because the fact that $|M/H|$
divides $|M_i/H|$ for two $i$ implies that $\alpha(M) \, | \, |H|$
as well.

Therefore, we assume that the sockel $S/H = M/H = A_k$.
Then as noticed by Huckleberry and Sauer, the conjugation map
induces an injection $1 \rightarrow G/H \rightarrow \Aut(S/H)$.
Hence we have the following, where the first
factor 2 is needed only when $k = 6$
(so that $\Aut(A_6)/A_6 = (\ZZ/(2))^{\oplus 2}$;
see Atlas \cite{Atlas}) :
$$|G/H| \le |\Aut(S/H)| \le 2 \times 2 \, |A_k| = 4 |M/H|.$$
First, treat the case $k \le 58$. Then $|G| \le 4|A_k| |H| \le
4|A_k| V \le 4(58!)V$ by Lemma \ref{Kquot}, and we are done. Or as
noticed by Szabo \cite{Sz}, $M/H$ contains a Sylow $p_1$-subgroup
$U/H$ such that (as above) $|M/H| \le |U/H|^5 \le (x_n
(V/|H|)^2)^5$, whence $|G| \le d_nV^{10}$ with $d_n = 4 x_n^5$.

We next deal with the case $k \ge 59$.
We use the approach of Szabo \cite{Sz}, but we apply the
centralizer Lemma \ref{gTh} instead.
As in Lemma \ref{alt}, set $\ell = [k/4]$
and we have $A_{\ell} \times A_{3\ell} \le A_k \le A_{4\ell + 3}$.

Suppose that the subgroup $A_{\ell} < A_k = M/H$ acts freely on
$X/H$. Then $\tilde{A}_{\ell} < G$ (with $\tilde{A}_{\ell}/H =
A_{\ell}$) also acts freely on $X$ by Lemma \ref{GGx}, whence
$|\tilde{A}_{\ell}| \le K_X^n = V$ by Lemma \ref{Kquot}. By Lemma
\ref{alt}, $|A_k| \le |A_{\ell}|^8 \le (V/|H|)^8$, so $|G| = |G/H|
|H| \le 4 |A_k| |H| \le 4V^8$.

Suppose that $A_{\ell}$ does not act freely on $X/H$.
We may take $\overline{\tau} \in G/H$ to be from $A_{\ell}$
and $\tau \in M$ with $\ord(\tau) = p = \ord(\overline{\tau})$.
Now $A_{3\ell} \le C_{M/H}(\overline{\tau}) \le C_{G/H}(\overline{\tau})$.
This and Lemma \ref{gTh} imply
$|A_{3\ell}| \le |C_G(\tau)| = |G_1|$.
By Lemmas \ref{alt}, \ref{Kquot} and \ref{central},
$|G| \le 4|A_k| |H| \le 4 |A_{3 \ell}|^{1.7} V
\le 4 |G_1|^{1.7} V \le 4(k_n q^2V^{3 + \varepsilon})^{1.7} V
= d_n q^{3.4} V^{6.1 + 1.7 \varepsilon}$,
where $d_n = 4 k_n^{1.7}$.
Since $p$ divides $|A_{\ell}|$, our $p^2$ divides
$A_k = M/H$ and also $|G|$.
So $\varepsilon = 0$ and we are done.
This completes the proof of Proposition \ref{key}.

\begin{setup}
{\bf Proof of Theorem \ref{Th2'}}
\end{setup}

Since $X$ is of general type, $G \le \Aut(X)$ is finite. Write
$\Aut_{\variety}(A) = T \rtimes A_0$ as at the beginning of the
section. Set $T_G = T \cap G$ which acts freely on $A$ and $X$. If
$G = T_G$, then $|G| \le K_X^n$ by Lemma \ref{Kquot}, and we are
done. So assume that $T_G < G$. Note that $A \rightarrow A/T_G$ is
an isogeny of abelian varieties, $X \rightarrow X/T_G$ is etale and
$G/T_G \le \{g \in \Aut_{\variety}(A/T_G) \, | \, g(X/T_G) =
X/T_G\}$. We shall check the conditions in \ref{setup1'} for $G/T_G
\le \Aut(X/T_G)$ and $X/T_G \subset A/T_G$. For every $g \in G
\setminus T_G$, we have $A^g \ne \emptyset$ by the assumpton on $G$,
so $(A/T_G)^{\overline g} \ne \emptyset$ by Lemma \ref{GGx}; here
$\overline{g} = g T_G \in G/T_G$. Thus $\overline{g}$ is not a
translation on $A/T_G$. So all conditions in \ref{setup1'} are
satisfied by $G/T_G$ and $X/T_G \subset A/T_G$; see Lemma \ref{GGx}
and the assumption on $G$. By Proposition \ref{key}, we have
$|G/T_G| \le d_n q(A/T_G)^{10-b} (K_{X/T_G}^n)^b = d_n q^{10-b}
(V/|T_G|)^b$. Now the theorem follows because $b \ge 5$. This proves
Theorem \ref{Th2'}.

\begin{setup}
{\bf Proof of Corollary \ref{Th2''}} \end{setup}

By Lemma \ref{sA}, $X$ is of general type and the fixed locus $A^g$
for every $\id \ne g \in G$ is a non-empty finite set unless $g$ is
a translation of $A$. Thus we can apply Theorem \ref{Th2'}. Now
Corollary \ref{Th2''} with $d_n' = d_n e_n^5$, follows from Theorem
\ref{Th2'}, Lemma \ref{sA}, Proposition \ref{qp} and Theorem
\ref{Gauss}.
\section{\bf Proofs of Theorems \ref{Th1} and \ref{Th2}}
In this section we shall prove Theorems \ref{Th1} and \ref{Th2}. We
start with an observation which bounds $\deg(\alb_X)$ in terms of
the volume of $X$. Note that $\Alb(X) = \Alb(X')$ if $X$ and $X'$
are smooth and birational. The lemma below follows from Lemma \ref{Kdeg} and
Theorem \ref{Gauss} (1) after resolving singularities of $W$,
since the volumes are birational invariants.

\begin{lem}\label{Kdeg'}
Let $X$ be a smooth projective $n$-fold. Suppose that $\alb_X : X
\rightarrow W := \alb(X) \subset \Alb(X)$ is generically finite onto
$W$, and $W$ is of general type. Then $V(X) \ge
\deg(\alb_X) \,\, V(W) \ge \deg(\alb_X)$.
\end{lem}

\begin{setup}
{\bf Proof of Theorem \ref{Th2}}
\end{setup}

Let $X$ and $G = \Bir(X)$ be as in Theorem \ref{Th2}. By Hanamura
\cite{Han} Lemma 2.4, after a smooth modification of $X$, we may
assume that $G$ acts regularly on $X$. The universal property of
$\Alb(X)$ implies that there is an induced action (not necessarily
faithful) of $G$ on $\Alb(X)$ such that $\alb_X : X \rightarrow
\Alb(X)$ is $G$-equivariant. Let $K$ (resp. $\overline{G}$) be the
kernel (resp. image) of the homomorphism $G = \Aut(X) \rightarrow
\Aut_{\variety}(\Alb(X))$ so that we have the exact sequence
$$1 \rightarrow K \rightarrow G
\rightarrow \overline{G} \rightarrow 1.$$ Since $K$ acts trivially
on $W = \alb_X(X)$, we can factor $\alb_X$ as $X \rightarrow X/K
\rightarrow W$. In particular, $|K| \le \deg(\alb_X)$. By Lemma
\ref{sA} and the assumption in Theorem \ref{Th2}, the fixed locus
$\Alb(X)^{\overline g}$ for every $\id \ne \overline{g} \in
\overline{G}$ is a non-empty finite set unless $\overline{g}$ is a
translation of $\Alb(X)$. Thus, by Theorem \ref{Gauss} (2), we can
apply Proposition \ref{qp} and Theorem \ref{Th2'} to $W$. Setting $V
= V(X)$ and $V(W) = K_W^n$ and noting that $q: = q(W) = q(X) = \dim
\Alb(X)$ by the universal property and definition of $\Alb(X)$, one
obtains the inequalities below with $d_n'' = d_n e_n^{10-b}$:
$$|\overline{G}| \le d_n q^{10-b} V(W)^b \le d_n'' V(W)^{10}.$$
Now by Lemma \ref{Kdeg'} or \ref{Kdeg}, we conclude Theorem \ref{Th2}:
$$|G| = |K| |\overline{G}| \le \deg(\alb_X) |\overline{G}|
\le d_n'' V \cdot V(W)^{9} \le d_n'' V^{10}.$$

\begin{setup}
{\bf Proof of Theorem \ref{Th1}} \end{setup}

Theorem \ref{Th1} (1) is a special case of Theorem \ref{Th1} (2).
Therefore, we have only to show Theorem \ref{Th1} (2).

So suppose that $\alb_X$ is not generically finite onto a 3-fold of
general type. Set $W = \alb_X(X) \subset \Alb(X)$. Then either $\dim
W < 3$, or $\kappa(W) < \dim W \le 3$. As above, after a smooth
modification of $X$, we may assume that $G = \Bir(X)$ acts regularly
on $X$ and of course on $\Alb(X)$ (not necessarily faithful on the
latter) so that $\alb_X : X \rightarrow \Alb(X)$ is $G$-equivariant.
Since $\dim \Alb(X) = q(X) \ge 4 > \dim X \ge \dim W$ by the
assumption, our $W$ is a proper subvariety of $\Alb(X)$. Hence by
Ueno \cite{Ue} Lemma 9.14 and Corollary 10.4, the Kodaira dimension
$\kappa(W) \ge 1$. As in Ueno \cite{Ue} Theorem 10.9 or Mori
\cite{Mo87} Theorem 3.7, let $B$ be the identity connected component
of $\{a \in \Alb(X) \, | \, a + W \subseteq W\}$. Then $W/B$ is of
general type, $W \rightarrow W/B$ is an etale fibre bundle with
fibre $B$ and is birational to the Iitaka fibring of $W$. Since the
pluri-canonical systems of $W$ are $G$-stable, we may take
$G$-equivariant desingularizations $X' \rightarrow X$ and $Y'
\rightarrow W/B$ such that $X' \rightarrow Y'$ is a well-defined
$G$-equivariant morphism, though $G$ might not act faithfully on
$Y'$. Rewrite $X'$ as $X$. Take a Stein factorization $X \rightarrow
Y \rightarrow Y'$ so that $f: X \rightarrow Y$ has connected general
fibre $F$ and is necessarily $G$-equivariant. We may also assume
that $Y$ is already smooth (or do equivariant modifications again).
Note that $G = \Bir(X) = \Aut(X)$ now. We have $\dim Y = \dim Y' =
\kappa(W) = 1, 2$. Since the subvariety $W/B < \Alb(X)/B$ is of
general type, both $Y$ and $Y'$ are of general type. Apply Theorem
\ref{Th3.1} to the $G$-equivariant map $f : X \rightarrow Y$ with
connected general fibre $F$ say. Note that $3 = \dim X = \dim Y +
\dim F$ and $|\Bir(Z)| = |\Aut(Z_{\min})| \le (42)^{\dim Z} V(Z)$
for both $Z = Y$ and $Z = F$, thanks to Hurwitz and Xiao. Here
$Z_{\min}$ is the {\it unique} smooth minimal model of $Z$, since
$\dim Z \le 2$. By Theorem \ref{Th3.1}, we have $|G| \le a \ V(X)$
with $a = (42)^3/3$. This proves Theorem \ref{Th1} (2).

\par \vskip 1pc \noindent
Department of Mathematics, National University of Singapore
\par \noindent
2 Science Drive 2, Singapore 117543, SINGAPORE
\par \noindent
E-mail: matzdq@nus.edu.sg

\newpage
\section{\bf Appendix}
$$\text{\bf A PRODUCT FORMULA FOR VOLUMES OF VARIETIES}$$

$$\text{By Yujiro Kawamata}$$

\par \vskip 1pc
The volume $v(X)$ of a smooth projective variety $X$ is defined by
\[
v(X) = \text{lim sup} \frac{\dim H^0(X, mK_X)}{m^d/d!}
\]
where $d = \dim X$.
This is a birational invariant.

\begin{thm}
Let $f: X \to Y$ be a surjective morphism of smooth projective
varieties with connected fibers. Assume that both $Y$ and the
general fiber $F$ of $f$ are varieties of general type. Then
\[
\frac{v(X)}{d_X!} \ge \frac{v(Y)}{d_Y!}\frac{v(F)}{d_F!}
\]
where $d_X = \dim X$, $d_Y = \dim Y$ and $d_F=\dim F$.
\end{thm}

\begin{proof}
Let $H$ be an ample divisor on $Y$.
There exists a positive integer $m_0$ such that $m_0K_Y - H$ is effective.

Let $\epsilon$ be a positive integer.
By Fujita's approximation theorem (\cite{F}),
after replacing a birational model of $X$,
there exists a positive integer $m_1$
and ample divisors $L$ on $F$ such that
$m_1K_F - L$ is effective and $v(\frac 1{m_1}L) > v(F)- \epsilon$.

By Viehweg's weak positivity theorem (\cite{V}),
there exists a positive integer $k$
such that $S^k(f_*\mathcal{O}_X(m_1K_{X/Y}) \otimes \mathcal{O}_Y(H))$ is
generically generated by global sections for a positive integer $k$.
$k$ is a function on $H$ and $m_1$.

We have
\[
\begin{split}
&\text{rank Im}(S^mS^k(f_*\mathcal{O}_X(m_1K_{X/Y})) \to
f_*\mathcal{O}_X(km_1mK_{X/Y})) \\
&\ge \dim H^0(F,kmL) \\
&\ge (v(F)- 2 \epsilon)\frac{(km_1m)^{d_F}}{d_F!}
\end{split}
\]
for sufficiently large $m$.

Then
\[
\begin{split}
&\dim H^0(X,km_1mK_X) \\
&\ge \dim H^0(Y,k(m_1-m_0)mK_Y) \times
(v(F)- 2 \epsilon)\frac{(km_1m)^{d_F}}{d_F!} \\
&\ge (v(Y) - \epsilon)\frac{(k(m_1-m_0)m)^{d_Y}}{d_Y!}
(v(F)- 2 \epsilon)\frac{(km_1m)^{d_F}}{d_F!} \\
&\ge (v(Y) - 2 \epsilon)(v(F)- 2 \epsilon)\frac{(km_1m)^{d_X}}{d_Y!d_F!}
\end{split}
\]
if we take $m_1$ large compared with $m_0$ such that
\[
\frac{(v(Y) - \epsilon)}{(v(Y)- 2 \epsilon)} \ge (\frac{m_1}{m_1-m_0})^{d_Y}.
\]
\end{proof}

\begin{rem}
If $X = Y \times F$, then we have an equality in the formula.
We expect that the equality implies the isotriviality of the family.
\end{rem}


\par \vskip 1pc \noindent
Department of Mathematical Sciences, University of Tokyo
\par \noindent
Komaba, Meguro, Tokyo, 153-8914, JAPAN
\par \noindent
E-mail: kawamata@ms.u-tokyo.ac.jp

\end{document}